\documentclass{elsart}

\usepackage{color}
\usepackage[T1]{fontenc}
\usepackage{subfigure}
\usepackage{graphicx}
\usepackage{amsmath,amssymb}
\usepackage{color}
\definecolor{oneblue}{rgb}{0,0.0,0.75}

\usepackage{dsfont}
\usepackage{textgreek}
\usepackage{mathrsfs}
\def\proof{{\noindent \it Proof. }}
\newtheorem{theorem}{Theorem}
\newtheorem{proposition}[theorem]{Proposition}
\newtheorem{lemma}[theorem]{Lemma}

\begin{document}
\begin{frontmatter}
\title{On a model for internal waves in rotating fluids}

\author{A. Dur\'an}
\address{ Applied Mathematics Department,  University of
Valladolid, P/ Bel\'en 15, 47011 Valladolid, Spain. Email:angel@mac.uva.es}


\begin{abstract}
In this paper a rotating two-fluid model for the propagation of internal waves is introduced. The model can be derived from a rotating-fluid problem by including gravity effects or from a nonrotating one by adding rotational forces in the dispersion balance. The physical regime of validation is discussed and mathematical properties of the new system, concerning well-posedness, conservation laws and existence of solitary-wave solutions, are analyzed.
\end{abstract}
\end{frontmatter}

\section{Introduction}
\label{sec1}
The purpose of this paper is to introduce and analyze a nonlinear, dispersive model for the one-way propagation of long internal waves of small amplitude along the interface of a two-layer system of fluids under gravity, surface tension and rotational effects. The model can be derived from different points of view; rotating- and nonrotating-fluid models proposed in the literature, mainly the Ostrovsky equation, \cite{Ostrovsky1978,Grimshaw1985}, and the Benjamin equation, \cite{ben0,ben1,ben2,abr}, respectively. The analysis of the new system, exposed in the present paper, is focused on mathematical aspects concerning well-posedness, conserved quantities and existence of solitary wave solutions. The main highlights are the following:
\begin{itemize}
\item Sufficient conditions on the parameters of the model are given in order to obtain existence and uniqueness of solutions of the associated linear problem. The result makes use of the theory on oscillatory integrals and regularity of dispersive equations developed in \cite{KenigPV1991b} (see also \cite{KenigPV1991,KenigPV1993,EsfahaniL2013,LevandoskyL2007,LiuV2004}).
\item The equation is shown to admit three conserved quantities by decaying to zero at infinity and smooth enough solutions. A Hamiltonian formulation is also derived.
\item One of the relevant properties of nonlinear dispersive models for wave propagation is the existence of traveling-wave solutions of solitary type, \cite{Bona1980,Bona1981} (see \cite{GilmanGS1995,Grimshaw1997} and references therein for the case of internal waves). In this sense, and using the Concentration-Compactness theory, \cite{Lions}, the new model is proved to admit such solutions, under suitable conditions on the parameters. By using the Petviashvili's iterative method, \cite{Pet1976}, to generate approximations to the solitary-wave profiles, several properties of the waves are analyzed by computational means. They concern the speed-amplitude relation, the asymptotic decay and the comparison with similar structures presented in classical rotating-fluid models like the Ostrovsky equation.
\end{itemize}
The paper is structured as follows. In Section~\ref{sec2}, the model will be introduced, from the general problem of propagation of internal waves along the interface of a two-layer system and under the corresponding physical regime of validation. Its justification from existing rotating- and nonrotating-fluid models by incorporating new physical assumptions will be discussed. 
Section~\ref{sec2} is finished off with the analysis of linear well-posedness of the corresponding initial-value problem (IVP) and the derivation of functionals preserved by smooth enough solutions vanishing suitably at infinity. In particular, a Hamiltonian structure of the problem comes out from one of these quantities. Section~\ref{sec3} is focused on the existence of solitary-wave solutions. As a first approach we make a computational study, with a description of the numerical technique used to generate approximate solitary-wave profiles and the numerical illustration of some of their properties. Then a theoretical result of existence of these solutions, under suitable conditions on the parameters of the model, is established.  These conditions for the existence will also help us to compare, by computational means, the proposed model with classical rotating-fluid models such as the Ostrovsky equation, with the aim of investigating the influence of the new physical properties assumed. Conclusions and future lines of research will be outlined in Section~\ref{sec4}.
%
%
%
%
%
%
%
%
%
%
%
%
%
%
%
%
%
%
%
%

The following notation will be used throughout the paper. By $H^{s}=H^{s}(\mathbb{R}), s\geq 0$ we denote the Sobolev space of order $s$, with $H^{0}=L^{2}(\mathbb{R})$ and norm denoted by $||\cdot ||_{s}$ (with $||\cdot ||=||\cdot||_{0}$). For $1\leq q\leq\infty, L^{q}=L^{q}(\mathbb{R})$ is the space of $q-$integrable functions with norm $||\cdot||_{L^{q}}$. On the other hand, $W^{1,q}(\mathbb{R}), q\geq 2$ (resp. $W_{loc}^{1,q}(\mathbb{R})$) will stand for the space of functions in $L^{q}$ (resp. locally in $L^{q}$) with weak derivative in $L^{q}$ (resp. locally in $L^{q}$).
\section{The mathematical model}
\label{sec2}
\subsection{On the derivation}

\begin{figure}[h!]
\centering
\includegraphics[width=0.8\textwidth]{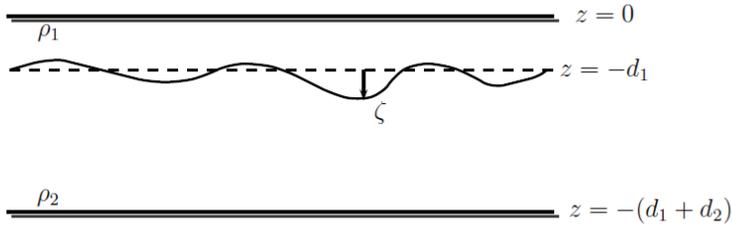}
\caption{Idealized model of internal wave propagation in a two-layer interface. $\rho_{2}>\rho_{1}; d_{2}>d_{1}$; $\zeta(x,t)$ denotes the downward vertical displacement of the interface from its level of rest at $(x,t)$.}
\label{amns_fig1}
\end{figure}

%
%
%
The two-layer interface problem for internal wave propagation, of interest for the present paper, is idealized in Figure~\ref{amns_fig1}.
This consists of two inviscid, homogeneous, incompressible fluids of depths $d_{i}, i=1,2$, with $d_{2}>d_{1}$ and densities $\rho_{i}, i=1,2$ with $\rho_{2}>\rho_{1}$. The upper and lower layers are respectively bounded above and below by a rigid horizontal plane, while the deviation of the interface from a level of rest, denoted by $\zeta$, is supposed to be a graph over the bottom.

\begin{figure}[h!]
\centering
\includegraphics[width=0.5\textwidth]{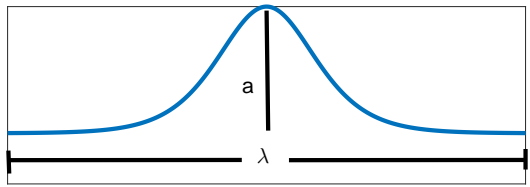}
\caption{Amplitude ($a)$ and wavelength ($\lambda$) of a wave.}
\label{amns_fig2}
\end{figure}

From this idealized system and in order to limit the physical regime of validation of the proposed model, some hypotheses are assumed. The first one is described in terms of the dimensionless parameters as
$$\quad \epsilon:=\frac{a}{d_{1}}<<1, \; \mu:=\left(\frac{d_{1}}{\lambda}\right)^{2}<<1,$$
referred to the upper fluid layer and where $a$ and $\lambda$ denote, respectively, typical amplitude and wavelength of a wave, see Figure~\ref{amns_fig2}. Thus we are assuming that the waves considered are long and of small amplitude with respect to the upper layer. 
The dispersive and nonlinear effects, governed by the parameters $\mu$ and $\epsilon$ respectively, are assumed to be balanced in the form
$$\mu\sim \epsilon^{2}.$$
Finally, capillary and gravity forces are assumed to be nonnegligible, as well as a dispersion effect due to the rotations of the fluids. These assumptions are translated to the following partial differential equation (PDE) for the evolution of the deviation of the interface
\begin{equation}\label{amns_e1}
\left(\zeta_t+\alpha \zeta_x+\zeta \zeta_x-\beta \mathcal{H}\zeta_{xx}-\delta \zeta_{xxx}\right)_{x}
=\gamma \zeta,
\end{equation}
where $\zeta=\zeta(x,t), x\in \mathbb{R}, t\geq 0$, and $\alpha, \gamma, \delta\geq 0,\; \beta\neq 0$ are constants.  If $\mathcal{H}$ denotes the Hilbert transform on $\mathbb{R}$,
\begin{equation}\label{amns_e2}
\mathcal{H}f(x):=\frac{1}{\pi}p.v.\int_{-\infty}^{\infty}\frac{f(y)}{x-y}\,dy,
\end{equation}
being one of the nonlocal terms of (\ref{amns_e1}), then the general dispersive effects (surface tension and gravity) are controlled by the parameters $\beta$ and $\delta$ (which depend on $\epsilon, \mu$, the densities $\rho_{i}, i=1,2$ as well as the interfacial surface tension and the acceleration of gravity, \cite{abr}). The dispersion due to the rotation is governed by the parameter $\gamma$ while $\alpha$ depends on the densities of the fluids. (One can always assume $\alpha=0$ by using the change of variables $x\mapsto x-\alpha t, t\mapsto t, \zeta\mapsto \zeta$.) The nonlinear effects are supposed to be of quadratic type.

The equation (\ref{amns_e1}) includes some well-known limiting cases from which its derivation can be justified. These can be rotating or nonrotating models. In the first case, one may start from the Ostrovsky equation, \cite{ApelOS2006,Grimshaw1985,GilmanGS1995,GilmanGS1996,
Grimshaw1997,Ostrovsky1978,
OstrovskyS1990,
Shira1981,Shira1986}
\begin{equation}\label{amns_e3}
\left(\zeta_t+\alpha \zeta_x+\zeta \zeta_{x}-\delta \zeta_{xxx}\right)_{x}=\gamma \zeta,
\end{equation}
and include the hypothesis of a much larger density of the lower fluid, $\rho_{2}>>\rho_{1}$, which implies a relevant gravity effect represented by $\beta$ and the Hilbert transform. Alternatively, one may consider the rotation-modified Benjamin-Ono (RMBO) equation, \cite{Grimshaw1985,GalkinS1991,LinaresM2004}
\begin{equation}\label{amns_e4}
\left(\zeta_t+\alpha \zeta_x+\zeta \zeta_{x}-\beta \mathcal{H}\zeta_{xx}\right)_{x}=\gamma \zeta,
\end{equation} 
and incorporate nonnegligible surface tension effects through the term associated to $\delta$. On the other hand, the most natural nonrotating model from which (\ref{amns_e1}) can be derived may be the Benjamin equation, \cite{ben0,ben1,ben2,abr}
\begin{equation}\label{amns_e5}
\zeta_t+\alpha \zeta_x+\zeta\zeta_x-\beta \mathcal{H}\zeta_{xx}-\delta \zeta_{xxx}
=0,
\end{equation}
if we assume (as in the case of the Ostrovsky equation with respect to the Korteweg-de Vries (KdV) equation, see \cite{Grimshaw1985}) that the rotational effects in the fluids are relevant enough to be included as a second nonlocal dispersive term $\gamma \partial_{x}^{-1}\zeta$, where $\partial_{x}^{-1}$ is defined as
\begin{equation*}
\partial_{x}^{-1}f(x)=\frac{1}{2}\left(\int_{-\infty}^{x}f(z)dz-\int_{x}^{\infty}f(z)dz\right),
\end{equation*}
or equivalently
\begin{equation}\label{amns_e6}
\widehat{\partial_{x}^{-1} f}(k)=(ik)^{-1}\widehat{f}(k),\quad k\in \mathbb{R}\backslash \{0\},\quad
\widehat{\partial_{x}^{-1} f}(0)=0,
\end{equation}
where
$$\widehat{f}(k):=\int_{-\infty}^{\infty} e^{-ikx}f(x)dx,\quad k\in \mathbb{R},
$$ denotes the Fourier transform of $f\in L^{2}$. Note that (\ref{amns_e6}) requires 
\begin{equation*}\label{amns_e7}
\widehat{f}(0):=\int_{-\infty}^{\infty} f(x)dx=0.
\end{equation*}
Due to the relation with the Benjamin equation (\ref{amns_e5}), equation (\ref{amns_e1}) will be sometimes referred as the rotation-modified Benjamin (RMBenjamin) equation.

We finally observe that, as in the case of other models, \cite{LevandoskyL2006}, an extension of (\ref{amns_e1}) is obtained by considering general homogeneous nonlinearities $f$ of some degree $p>1$, that is
\begin{equation}\label{rmben0}
f(\lambda s)=\lambda^{p}f(s),\quad s, \lambda\neq 0.
\end{equation}
 in such a way that (\ref{amns_e1}) can be generalized to
\begin{equation}\label{rmben1}
(u_{t}+\alpha u_{x}+f(u)_{x}-\beta\mathcal{H}u_{xx}-\delta u_{xxx})_{x}=\gamma u,\quad x\in\mathbb{R}, \quad t>0.
\end{equation}
The main theoretical results below will be established for (\ref{rmben1}), although the particular case of (\ref{amns_e1}) (for which $f(u)=u^{2}/2$ and $p=2)$) may be of more interest. 

\subsection{Well-posedness}
This section concerns the well-posedness of the IVP of the linearized equation associated to (\ref{rmben1})
\begin{align}
\left(\zeta_t+\alpha \zeta_x-\beta \mathcal{H}\zeta_{xx}-\delta \zeta_{xxx}\right)_{x}=\gamma \zeta,\label{amns_e8}\\
\zeta(x,0)=\zeta_{0}.\nonumber
\end{align}
Using the Fourier representation of (\ref{amns_e2})
\begin{equation*}
\widehat{\mathcal{H}f}(k)=-i{\rm sign}(k)\widehat{f}(k), \quad k\in\mathbb{R},
\end{equation*}
then the application of the Fourier transform (in $x$) to (\ref{amns_e8}) leads to
\begin{equation}\label{amns_e9}
\widehat{\zeta}(k,t)=e^{-im(k)t}\widehat{\zeta_{0}}(k),\quad k\in\mathbb{R}
\end{equation}
where
\begin{equation*}\label{amns_e10a}
m(k)=\frac{\gamma}{k}+\alpha k-\beta k|k|+\delta k^{3}, k\neq 0,\quad m(0)=0.
\end{equation*}
The inversion of (\ref{amns_e9}) allows to write formally the solution of (\ref{amns_e8}) in the operational form
\begin{equation}\label{amns_e11a}
\zeta(x,t)=S(t)\zeta_{0}(x)=\frac{1}{2\pi}\int_{-\infty}^{\infty} e^{i(kx-m(k)t)}\widehat{\zeta_{0}}(k)dk.
\end{equation}
The following lemma will be used to estimate (\ref{amns_e11a}).
\begin{lemma}
\label{lemma1}
For $\gamma, \delta>0$, $\beta$ satisfying $\beta<0$ or $0<\beta <4\gamma^{1/4}\delta^{3/4}$ and
\begin{equation}\label{amns_e12a}
\phi(k)=\frac{\gamma}{k}-\beta k|k|+\delta k^{3}, k\neq 0,\quad \phi(0)=0,
\end{equation}
then $|\phi^{\prime\prime}(k)|\geq -2\beta+8\gamma^{1/4}\delta^{3/4}, k\neq 0$.
\end{lemma}
\proof Let $k_{0}=\left(\gamma/\delta\right)^{1/4}$. One can check that if $k<0$ then
$$\phi^{\prime\prime}(k)\geq \phi^{\prime\prime}(-k_{0})=-2\beta+8\gamma^{1/4}\delta^{3/4},$$
while if $k>0$ then
$$\phi^{\prime\prime}(k)\leq \phi^{\prime\prime}(k_{0})=2\beta-8\gamma^{1/4}\delta^{3/4}.$$ These two inequalities, under the hypotheses on $\beta$, prove the result.$\Box$

Following similar arguments to those of \cite{LinaresM2004,LinaresM2006,VarlamovL2004} for the case of the Ostrovsky (\ref{amns_e3}) and the RMBO equations,  (\ref{amns_e3}) and (\ref{amns_e4}) respectively, we have the following result on well-posedness of (\ref{amns_e8}).
\begin{theorem}
\label{theorem1}
Under the hypotheses of Lemma~\ref{lemma1}, let $f\in L^{2}$ and $t>0$. Then
\begin{equation}\label{amns_e13aa}
||S(t)f||_{L_{t}^{q}L_{x}^{p}}\leq C||f||_{L^{2}},
\end{equation}
for some $C>0$ and where $q=\frac{4}{\theta}, p=\frac{2}{1-\theta}, \theta\in [0,1]$.
\end{theorem}

\proof Recall that a change of variable allows to assume $\alpha=0$ in (\ref{amns_e8}). Note that under the hypotheses on $\beta, \gamma, \delta$ of Lemma~\ref{lemma1}, the function $\phi$ in (\ref{amns_e12a}) satisfies the conditions (2.1a)-(2.1e) described in \cite{KenigPV1991b} and applying Theorem 2.1 of this reference, we have
\begin{equation*}
||W_{\theta/2}(t)f||_{L_{t}^{q}(\mathbb{R},L^{p})}\leq C||f||_{L^{2}},
\end{equation*}
for some constant $C$ and where
\begin{equation*}
W_{s}(t)f(x):=\int_{-\Omega} e^{i(kx+\phi(k)t)}|\phi^{\prime\prime}(k)|^{s/2}\widehat{f}(k)dk, \quad s\geq 0,
\end{equation*}
with $\Omega=(-\infty,0)\cup (0,\infty)$. Now, using Lemma~\ref{lemma1}, observe that (cf. \cite{LinaresM2004,LinaresM2006})
\begin{equation*}
||S(t)f||_{L_{t}^{q}L_{x}^{p}}\leq \frac{c_{0}}{2\pi}||W_{\theta/2}(t)f||_{L_{t}^{q}(\mathbb{R},L^{p})},\quad c_{0}=\frac{1}{-2\beta+8\gamma^{1/4}\delta^{3/4}},
\end{equation*}
and (\ref{amns_e13aa}) follows.$\Box$
\begin{note}
It may be worth emphasizing the particular cases of (\ref{amns_e13aa}) corresponding to the limiting values $\theta=0, 1$:
\begin{equation*}
||S(t)f||_{L_{t}^{\infty}L_{x}^{2}}\leq C||f||_{L^{2}},
\quad
||S(t)f||_{L_{t}^{2}L_{x}^{\infty}}\leq C||f||_{L^{2}}.
\end{equation*}
\end{note}
\subsection{Conserved quantities}
A second mathematical property is concerned with the existence of invariant quantites of (\ref{rmben1}) for smooth enough solutions. Similar results to those of other rotating models, like the Ostrovsky equation, or nonrotating fluid models like the Benjamin equation, can be derived in this case. The proof is straightforward.
\begin{proposition}
\label{proposition1}
Assume that $\zeta$ is a smooth solution of (\ref{rmben1}) such that $$\zeta,\zeta_{x},\zeta_{xx},\zeta_{xxx},\zeta_{xxxx}\rightarrow 0,$$ as $x\rightarrow\pm\infty$. Then $\zeta$ satisfies the zero mass condition
\begin{equation}\label{amns_e14a}
I(\zeta)=\int_{-\infty}^{\infty}\zeta(x,t)dx=0,\quad t\geq 0
\end{equation}
and the time preservation of the momentum and energy
\begin{eqnarray}
V(\zeta)&=&\int_{-\infty}^{\infty}\frac{\zeta(x,t)^{2}}{2}dx,\nonumber\\
E(\zeta)&=&\int_{-\infty}^{\infty}\left(\alpha\frac{\zeta(x,t)^{2}}{2}+F(\zeta(x,t))-\frac{\beta}{2}\zeta(x,t)\mathcal{H}\zeta_{x}(x,t)\right.\nonumber\\
&&\left.+\frac{\delta}{2}(\zeta_{x}(x,t))^{2}+\frac{\gamma}{2}((\partial_{x}^{-1}\zeta)(x,t))^{2}\right)dx, \label{rmben2}
\end{eqnarray}
where $F^{\prime}=f, F(0)=0$. The functional (\ref{rmben2}) is the Hamiltonian of (\ref{rmben1}) with respect to the symplectic structure given by $J=-\partial/\partial x$.
\end{proposition}

 \section{Solitary wave solutions}
\label{sec3}
A third mathematical property of (\ref{rmben1}), under study in the present paper, is concerned with the existence of solitary wave solutions. These are solutions of permanent form $\zeta=\varphi(x-c_{s}t)$ that travel with constant speed of propagation $c_{s}\neq 0$ and decay, along with its derivatives, to zero as $X=x-c_{s}t\rightarrow\pm \infty$. Substituting into (\ref{rmben1}) and integrating once, the profile $\varphi=\varphi(X)$ must satisfy
\begin{eqnarray}
(-c_{s}+\alpha)\varphi+f(\varphi)-\beta \mathcal{H}\varphi^{\prime}-\delta\varphi^{\prime\prime}
-\gamma\partial_{x}^{-2}\varphi=0,\label{rmben3}
\end{eqnarray}
where $\partial_{x}^{-2}\varphi:=\partial_{x}^{-1}(\partial_{x}^{-1}\varphi)$
\subsection{Numerical generation}
The numerical approximation to (\ref{rmben3}) may give us a first approach about the existence of solitary wave solutions and some of their properties. In this section this will be illustrated for the case of (\ref{amns_e1}), that is when $f(u)=u^{2}/2$ in (\ref{rmben1}). To this end, a typical strategy consists of considering the profiles as solutions of the fixed point equation
\begin{equation}\label{amns_e311}
\underbrace{(c_{s}-\alpha)\varphi+\beta \mathcal{H}\varphi^{\prime}+\delta\varphi^{\prime\prime}
+\gamma\partial_{x}^{-2}\varphi}_{\mathcal{L}\varphi}=\underbrace{\frac{\varphi^{2}}{2}}_{\mathcal{N}(\varphi)},
\end{equation}
which may be solved iteratively. Among other alternatives presented in the literature, (see e.~g. the review and references in \cite{yang2}), the numerical resolution of (\ref{amns_e311}) will be here performed with the Petviahsvili's method, \cite{Pet1976}. This is formulated as follows. Given an initial profile $\varphi^{[0]}$, the approximation $\varphi^{[\nu+1]}$ is computed from $\varphi^{[\nu]}, \nu=0,1,\ldots$, by
\begin{eqnarray}
m^{[\nu]}&=&\frac{\langle \mathcal{L}\varphi^{[\nu]},\varphi^{[\nu]}\rangle}{\langle \mathcal{N}(\varphi^{[\nu]}),\varphi^{[\nu]}\rangle},\label{amns_e312a}\\
\mathcal{L}\varphi^{[\nu+1]},\varphi^{[\nu]}&=&\left(m^{[\nu]}\right)^{2}\mathcal{N}(\varphi^{[\nu]}),\quad \nu=0,1,\ldots, \label{amns_e312b}
\end{eqnarray}
where $\mathcal{L}, \mathcal{N}$ are, respectively, the linear and nonlinear operators defined in (\ref{amns_e311}) and $\langle\cdot,\cdot,\rangle$ denotes the Euclidean inner product. The application of the Petviashvili's method requires, among other conditions, a nonlinearity $\mathcal{N}$ of homogeneous type. Its degree of homogeneity determines the exponent of $m^{[\nu]}$ in (\ref{amns_e312b}). This term is usually called the stabilizing factor and governs the convergence of the iteration, see e.~g. \cite{PelinovskyS2004}. As for the implementation of (\ref{amns_e312a}), (\ref{amns_e312b}), the equation (\ref{amns_e311}) is typically discretized on a long enough interval $(-L,L)$ by a Fourier pseudospectral approximation to the values of the profiles at a uniform grid $\{x_{j}=-L+jh, j=0,\ldots,N\}$ of collocation points, where $N\geq 1$ is an integer and $h:=2L/N$. The vector $\varphi_{h}=(\varphi_{h,0},\ldots,\varphi_{h,N-1})^{T}$ where $\varphi_{h,j}$ approximates $\varphi(x_{j}), j=0,\ldots,N-1$, satisfies a nonlinear algebraic system
\begin{equation}\label{amns_e313}
{\mathcal{L}}_{h}\varphi_{h}={\mathcal{N}}_{h}(\varphi_{h}),
\end{equation}
which is obviously an approximation to (\ref{amns_e311}) with approximate operators ${\mathcal{L}}_{h}, {\mathcal{N}}_{h}$ to ${\mathcal{L}}, {\mathcal{N}}$ respectively. The system (\ref{amns_e313}) is typically solved in Fourier space by using the Fourier symbol of the linear part, in such a way that the equations for the discrete Fourier coefficients of the approximation $\varphi_{h}$ will have the form
\begin{equation}\label{amns_e314}
p(k)\widehat{\varphi_{h}}(k)=k^{2}\widehat{\mathcal{N}_{h}(\varphi_{h})}(k),\; k=-N/2,\ldots,N/2-1, k\neq 0,\quad \widehat{\varphi_{h}}(0)=0,
\end{equation}
where in (\ref{amns_e314}), $\widehat{\varphi_{h}}(k)$ is the $k$-th discrete Fourier coefficient of $\varphi_{h}$, $p(k):=k^{2}(c_{s}-\alpha)+\beta k^{3}|k|-\delta k^{4}-\gamma$ and we choose ${\mathcal{N}}_{h}(\varphi_{h})=\frac{1}{2}\varphi_{h}.^{2}$, with the dot standing for the Hadamard product (componentwise) of vectors. The requirement $\widehat{\varphi_{h}}(0)=0$ is nothing but the zero mass condition (\ref{amns_e14a}). Thus, the Fourier representation of the discrete version of the Petviashvili's method (\ref{amns_e312a}), (\ref{amns_e312b}) has the following form
\begin{eqnarray}
m_{h}^{[\nu]}&=&\frac{\displaystyle\sum_{k}p(k)|\widehat{\varphi_{h}^{[\nu]}}(k)|^{2}}{ \displaystyle\sum_{k}\frac{k^{2}}{2}\widehat{(\varphi_{h}^{[\nu]})^{2}}(k)\overline{\widehat{\varphi_{h}^{[\nu]}}(k)}},\label{amns_e315a}\\
\widehat{\varphi_{h}^{[\nu+1]}}(k)&=&\left(m_{h}^{[\nu]}\right)^{2}\frac{k^{2}}{2p(k)}\widehat{(\varphi_{h}^{[\nu]})^{2}}(k),\; k\neq 0,\; \nu=0,1,\ldots,\label{amns_e315b}
\end{eqnarray}
with $\widehat{\varphi_{h}^{[\nu]}}(0)=0, \nu=0,1,\ldots,$ and from an initial iteration $\varphi_{h}^{[0]}$.

The purpose of the implementation of (\ref{amns_e315a}), (\ref{amns_e315b}) is two-fold: having more certainty about the existence of solitary waves and deriving a computational way to obtain approximate profiles from which some properties of the waves and their dynamics can be discussed. 
\begin{figure}[h!]
\centering
{\includegraphics[width=0.8\textwidth]{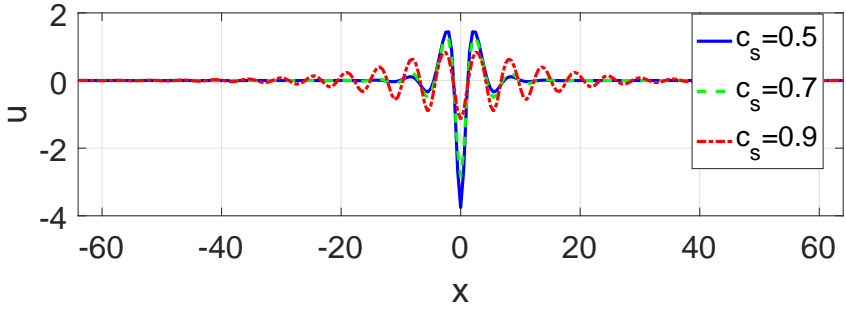}}
\caption{Numerical approximation with $\alpha=0, \beta=\gamma=\delta=1$. Computed solitary wave profiles.}
\label{amns_fig3}
\end{figure}
Thus, Figure~\ref{amns_fig3} shows some computed profiles corresponding to different values of the speed and for $\alpha=0, \beta=\gamma=\delta=1$. Two properties are suggested: the amplitude of the waves is a decreasing function of the speed and the waves have an oscillatory decay, with the oscillations increasing with the speed. 
\begin{figure}[h!]
\centering
{\includegraphics[width=0.8\textwidth]{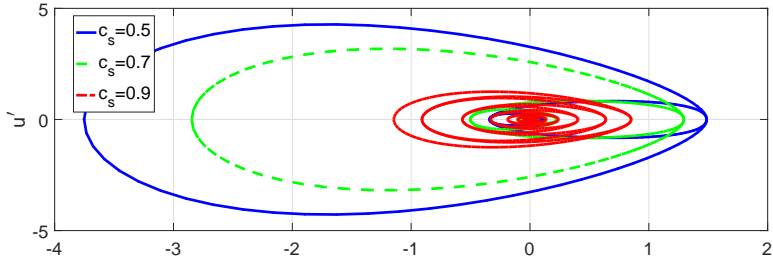}}
\caption{Numerical approximation with $\alpha=0, \beta=\gamma=\delta=1$. Phase portraits of the computed solitary wave profiles.}
\label{amns_fig4}
\end{figure}

These two properties are confirmed by the following figures. Figure~\ref{amns_fig4} displays the phase portraits of the profiles computed in Figure~\ref{amns_fig3} and the oscillatory decay is clearly observed. By fitting the values close to the origin, the results suggest that the waves decay algebraically, as in the cases of the Benjamin equation, \cite{abr} and the Ostrovsky equation, \cite{GalkinS1991,GilmanGS1995}.

\begin{figure}[h!]
\centering
{\includegraphics[width=0.8\textwidth]{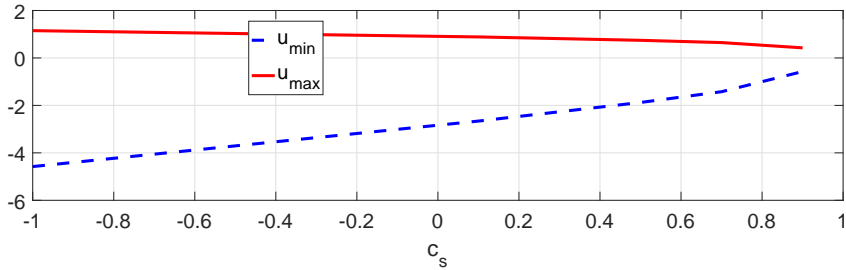}}
\caption{Numerical approximation with $\alpha=0, \beta=\gamma=\delta=1$. Speed-amplitude relations.}
\label{amns_fig5}
\end{figure}

Figure~\ref{amns_fig5} shows the behaviour of the maximum positive excursion $u_{max}$ and the minimum negative excursion $u_{min}$ of the profiles as functions of the speed. This confirms how the amplitude of the waves decreases as the speed is increasing. Additionally, Figure~\ref{amns_fig5} also suggests a limiting value of the speed $c_{s}$, far from which the existence of the solitary-wave profiles does not seem to be guaranteed.


\begin{note}
in the case $\beta<0$, the method (\ref{amns_e315a}), (\ref{amns_e315b}) can also compute multi-pulses. One of these is shown in Figure~\ref{amns_fig6}.
\begin{figure}[h!]
\centering
{\includegraphics[width=0.7\textwidth]{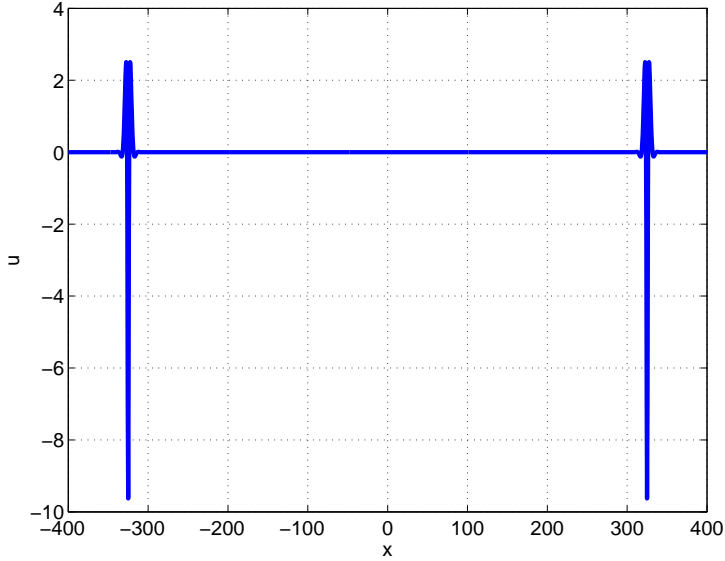}}
\caption{Two-pulse for $\alpha=0, \beta=-1,\gamma=\delta=1, p=2, c_{s}=1.1$ and a negative hyperbolic-secant profile as initial data for the iteration (\ref{amns_e315a}), (\ref{amns_e315b}).}
\label{amns_fig6}
\end{figure}
\end{note}
These computations motivate to study the existence of solitary wave solutions of (\ref{rmben1}) theoretically. This is developed in the following subsection.
\subsection{Existence of solitary wave solutions}
We consider the space
\begin{eqnarray*}
H=\left\{\varphi\in H^{3/2}(\mathbb{R})/ \partial_{x}^{-1}\varphi\in L^{2}(\mathbb{R})\right\},
\end{eqnarray*}
(in \cite{LiuV2004} this is denoted by $X_{3/2}$) where $\partial_{x}^{-1}$ is given by (\ref{amns_e6}). In $H$, the norm is defined as (cf. \cite{LiuV2004,EsfahaniL2013})
\begin{eqnarray*}
||\varphi||_{H}=||\varphi||_{H^{3/2}}+||\partial_{x}^{-1}\varphi||,
\end{eqnarray*}
or, equivalently
\begin{eqnarray*}
||\varphi||_{H}=||\partial_{x}\varphi||+||D_{x}^{1/2}\varphi||+||\partial_{x}^{-1}\varphi||,
\end{eqnarray*}
(see (\ref{amns_e11})), where $D_{x}^{1/2}$ is defined from the corresponding Fourier symbol $v(\xi)=|\xi|^{1/2},\xi\in\mathbb{R}$.

The purpose in this section is to discuss the existence of solitary-wave solutions of (\ref{rmben1}) in terms of the parameters $\alpha, \beta, \gamma$ and $\delta$. The discussion is based on the corresponding studies in the literature for the Ostrovsky equation, \cite{LiuV2004,LevandoskyL2006}, and generalized versions, \cite{LevandoskyL2007}, as well as the RMBenjamin-Ono equation, \cite{EsfahaniL2013}. In all the cases, some results of existence are obtained by applying the Concentration-Compactness theory, \cite{Lions}.

We define the functionals
\begin{eqnarray}
I(u)&=&\int_{-\infty}^{\infty}\left(-(c_{s}-\alpha)u^{2}-\beta u\mathcal{H}u_{x}+\gamma(\partial_{x}^{-1}u)^{2}+\delta u_{x}^{2}\right)dx;\label{rmben4a}\\
K(u)&=&-(p+1)\int_{-\infty}^{\infty}F(u)dx,\label{rmben4b}
\end{eqnarray}
and consider, for $\lambda>0$, the minimization problem
\begin{eqnarray}
M_{\lambda}=\inf\{I(u): u\in H, K(u)=\lambda\}.\label{rmben5}
\end{eqnarray}

\begin{note}
Note that if $\psi\in H$ achieves the minimum (\ref{rmben5}) for some $\lambda>0$, then there is a Lagrange multiplier $\mu\in\mathbb{R}$ such that $I(\psi)=\mu K(\psi)$. This means that
\begin{eqnarray*}
-\beta \mathcal{H}\psi_{x}-(c_{s}-\alpha)\psi-\delta\psi_{xx}-\gamma\partial_{x}^{-2}\psi=-\mu (p+1)f(\psi).
\end{eqnarray*}
If $\mu^{\prime}=(p+1)\mu$ then $\varphi=(\mu^{\prime})^{\frac{1}{p-1}}\psi$ satisfies (\ref{rmben3}).
\end{note}
\begin{note}
\label{rem2}
Note that if we integrate (\ref{rmben0}), with respect to $\lambda$, between $0$ and $t$ we have
\begin{eqnarray*}
\frac{1}{s}F(t s)=\frac{t^{p+1}}{p+1}F^{\prime}(s),
\end{eqnarray*}
which, evaluated at $t=1$, implies that
\begin{eqnarray}
(p+1)F(s)=sF^{\prime}(s),\label{rmben5b}
\end{eqnarray}
and, therefore, $F$ is homogeneous of degree $p+1$. Some consequences of this are:
\begin{itemize}
\item The functional $K$ in (\ref{rmben4b}) is homogeneous of degree $p+1$. (Note that $I$ in (\ref{rmben4a}) is homogeneous of degree two.)
\item There exists $C>0$ such that
\begin{eqnarray}
|F(u)|\leq C |u|^{p+1}.\label{rmben6}
\end{eqnarray}
\end{itemize}
They will be used elsewhere.
\end{note}
\begin{note}
\label{remark32b}
We denote by $G=G(\alpha,\beta,\gamma,\delta,c_{s})$ the set of solutions of (\ref{rmben3}). From the homogeneity of $I$ and $K$, $u\in G$ also achieves the minimum 
\begin{eqnarray*}
m=m(\alpha,\beta,\gamma,\delta,c_{s})=\inf\{\frac{I(u)}{K(u)^{\frac{2}{p+1}}}; u\in H, K(u)>0\},
\end{eqnarray*}
and therefore $M_{\lambda}=\lambda^{\frac{2}{p+1}}m$. If we multiply (\ref{rmben1}) by $\phi$, use (\ref{rmben5b}) and integrate, we have $I(\varphi)=K(\varphi)$, in such a way that
\begin{eqnarray*}
G=\left\{\varphi\in H / I(\varphi)=K(\varphi)=m^{\frac{p+1}{p-1}}\right\}
\end{eqnarray*}
\end{note}
\begin{note}
Throughout the rest of the paper the following estimates will be used:
\begin{eqnarray}
||D_{x}^{1/2}u||^{2}\leq \epsilon^{2}||u||^{2}+\frac{1}{4\epsilon^{2}}||\partial_{x}u||^{2},\label{amns_e10}
\end{eqnarray}
for any $\epsilon>0$.
\begin{eqnarray}
||u||^{2}&\leq &C||\partial_{x}^{-1}u||^{1/2}||\partial_{x}u||^{1/2},\label{amns_e11}\\
||\partial_{x}^{-1}u||_{\infty}&\leq &C||\partial_{x}^{-1}u||^{1/2}||u||^{1/2},\nonumber
\end{eqnarray}
for some constant $C$, see e.~g. \cite{LiuV2004}.
\end{note}
In order to prove that $G$ is not empty, we need several previous results. The first one is given in the following lemma.
\begin{lemma}
\label{lemma2}
Assume that $\delta,\gamma>0$ and that one the following conditions holds:
\begin{itemize}
\item[(i)] $\beta<0, c_{s}-\alpha<0$.
\item[(ii)] $\beta<0, 0<c_{s}-\alpha<c^{*}=2\sqrt{\gamma\delta}$.
\item[(iii)] $\beta>0, c_{s}-\alpha\leq-\frac{\beta^{2}}{4\delta}$.
\item[(iv)] $\beta>0$ with $4\delta-\beta>0, \beta^{3}<\gamma (4\delta-\beta)^{2}$ and $0<c_{s}-\alpha<z_{+}$ where
$$z_{+}=\frac{1}{2}\left(-\beta\left(1+\frac{\beta}{4\delta}\right)+(4\delta-\beta)\sqrt{\frac{\gamma}{\delta}+\left(\frac{\beta}{4\delta}\right)^{2}}\right).$$
\end{itemize}
Then $M_{\lambda}>0$ for $\lambda>0$.
\end{lemma}

\proof 

The proof of Lemma \ref{lemma2} is based on the following estimates of $I(u)$:
\begin{itemize}
\item In the case of (i):
\begin{eqnarray}
I(u)&\geq &(-\beta)\int |D_{x}^{1/2}u|^{2}dx+\gamma\int |\partial_{x}^{-1}u|^{2}dx+\delta\int |\partial_{x}u|^{2}dx.\label{amns_e12}
\end{eqnarray}
\item In the case of (ii):
\begin{eqnarray}
I(u)&\geq &(-\beta)\int |D_{x}^{1/2}u|^{2}dx+(\gamma-\left(\frac{c_{s}-\alpha}{4\epsilon^{2}}\right) )\int |\partial_{x}^{-1}u|^{2}dx\nonumber\\
&&+(\delta-\epsilon^{2}(c_{s}-\alpha))\int |\partial_{x}u|^{2}dx.\label{amns_e13}
\end{eqnarray}
for some $\epsilon^{2}\in\displaystyle\left(\frac{c_{s}-\alpha}{4\gamma},\frac{\delta}{c_{s}-\alpha}\right)$.
\item In the case of (iii):
\begin{eqnarray}
I(u)&\geq &(\alpha-c_{s}-\beta\epsilon^{2}))\int |u|^{2}dx+\gamma\int |\partial_{x}^{-1}u|^{2}dx\nonumber\\
&&+\left(\delta-\frac{\beta}{4\epsilon^{2}}\right)\int |\partial_{x}u|^{2}dx,\label{amns_e14}
\end{eqnarray}
for some $\epsilon^{2}\in\displaystyle\left(\frac{\beta}{4\delta},\frac{\alpha-c_{s}}{\beta}\right)$.
\item In the case of (iv):
\begin{eqnarray}
I(u)&\geq &(\gamma-\epsilon^{2}(c_{s}-\alpha+\beta\epsilon^{2}))\int |\partial_{x}^{-1}u|^{2}dx\nonumber\\
&&+\left(\delta-\frac{\beta}{4\epsilon^{2}}-\left(\frac{c_{s}-\alpha+\beta\epsilon^{2}}{4\epsilon^{2}}\right)\right)\int |\partial_{x}u|^{2}dx,\label{amns_e15}
\end{eqnarray}
for some $\epsilon$ such that
\begin{eqnarray*}
&&\frac{\beta+c_{s}-\alpha}{4\delta-\beta}<\epsilon^{2}<X_{+}\\
&&X_{+}=\frac{1}{2}\left(-\left(\frac{c_{s}-\alpha}{\beta}\right)+\sqrt{\left(\frac{c_{s}-\alpha}{\beta}\right)^{2}+\frac{4\gamma}{\beta}}\right).
\end{eqnarray*}
\end{itemize}

The proof of (\ref{amns_e12})-(\ref{amns_e15}) is as follows.
\begin{itemize}
\item The proof of (\ref{amns_e12}) is trivial since $c_{s}-\alpha<0$.
\item For the proof of (\ref{amns_e13}) we write, \cite{LiuV2004}
\begin{eqnarray*}
\int u^{2}dx=-\int(\partial_{x}^{-1}u)(\partial_{x}u)dx.
\end{eqnarray*}
Then, for any $\epsilon^{2}>0$
\begin{eqnarray*}
-(c_{s}-\alpha)\int u^{2}dx&\geq &(c_{s}-\alpha)\left(-\epsilon^{2}\int |\partial_{x} u|^{2}dx-\frac{1}{4\epsilon^{2}}\int |\partial_{x}^{-1}u|^{2}dx\right).
\end{eqnarray*}
This is applied to $I(u)$, leading to (\ref{amns_e13}). Since $0<c_{s}-\alpha<2\sqrt{\gamma\delta}$ then $$\displaystyle\frac{c_{s}-\alpha}{4\gamma}<\displaystyle\frac{\delta}{c_{s}-\alpha},$$ and choosing $\epsilon^{2}\in\left(\displaystyle\frac{c_{s}-\alpha}{4\gamma},\displaystyle\frac{\delta}{c_{s}-\alpha}\right)$ ensures that all the terms in the right hand side of (\ref{amns_e13}) are positive.
\item For the proof of (\ref{amns_e14}), we use (\ref{amns_e10}) and similar arguments to those of the previous proof to choose $\epsilon^{2}$.
\item Proof of (\ref{amns_e15}): In this case, the same strategy as above is applied twice. First, we have
\begin{eqnarray*}
I(u)&\geq &-(c_{s}-\alpha+\beta\epsilon_{1}^{2})\int u^{2}dx+\left(\delta-\frac{\beta}{4\epsilon_{1}^{2}}\right)\int |\partial_{x} u|^{2}dx+\gamma\int|\partial_{x}^{-1}u|^{2}dx\\
&\geq &(c_{s}-\alpha+\beta\epsilon_{1}^{2})\left(-\epsilon_{2}^{2}\int |\partial_{x}^{-1}u|^{2}dx-\frac{1}{4\epsilon_{2}^{2}}\int |\partial_{x} u|^{2}dx\right)\\
&&+\left(\delta-\frac{\beta}{4\epsilon_{1}^{2}}\right)\int |\partial_{x} u|^{2}dx+\gamma\int|\partial_{x}^{-1}u|^{2}dx\\
&=&\left(\delta-\frac{\beta}{4\epsilon_{1}^{2}}-\left(\frac{c_{s}-\alpha+\beta\epsilon_{1}^{2}}{4\epsilon_{2}^{2}}\right)\right)\int |\partial_{x} u|^{2}dx\\
&&+\left(\gamma-\epsilon_{2}^{2}\left(c_{s}-\alpha+\beta\epsilon_{1}^{2}\right)\right)\int|\partial_{x}^{-1}u|^{2}dx.
\end{eqnarray*}
Now the two coefficients are positive when
\begin{eqnarray}
&&\gamma-\epsilon_{2}^{2}(c_{s}-\alpha)-\beta\epsilon_{1}^{2}\epsilon_{2}^{2}>0,\label{amns_e13a}\\
&&(4\delta-\beta)\epsilon_{1}^{2}\epsilon_{2}^{2}-
\beta\epsilon_{2}^{2}-(c_{s}-\alpha)\epsilon_{1}^{2}>0.\label{amns_e13b}
\end{eqnarray}
Note that (\ref{amns_e13b}) implies that we need $4\delta-\beta>0$. If we simplify by setting $\epsilon_{1}^{2}=\epsilon_{2}^{2}=\epsilon^{2}$ then the satisfaction of (\ref{amns_e13a}), (\ref{amns_e13b}) requires to choose $\epsilon^{2}$ within the range specified in (\ref{amns_e15}).
\end{itemize}
Once (\ref{amns_e12})-(\ref{amns_e15}) is proved, we can use (\ref{rmben6}), (\ref{amns_e10}) and the estimate (14) of \cite{EsfahaniL2013} to have
\begin{eqnarray}
\lambda=K(u)\leq C\left(||u||^{2}+||\partial_{x}u||^{2}+||\partial_{x}^{-1}u||^{2}\right)^{\frac{p+1}{2}},\label{amns_e16}
\end{eqnarray}
for some constant $C$. From (\ref{amns_e16}) and (\ref{amns_e11}) we have
\begin{eqnarray}
\lambda=K(u)\leq C\left(||\partial_{x}u||^{2}+||\partial_{x}^{-1}u||^{2}\right)^{\frac{p+1}{2}},\label{amns_e17}
\end{eqnarray}
for some constant $C$. Now, using (\ref{amns_e10}) and (\ref{amns_e11}) if necessary, we have that in all the cases (i) to (iv) the right hand side of (\ref{amns_e17}) can be bounded  by the right hand side of the corresponding estimate (\ref{amns_e12}) to (\ref{amns_e15}), in such a way that there exists $C=C(c_{s},p,\alpha,\beta,\gamma,\delta,\epsilon^{2})>0$ such that
\begin{eqnarray*}
\lambda=K(u)\leq (I(u))^{\frac{p+1}{2}},
\end{eqnarray*}
which implies
\begin{eqnarray*}
I(u)\geq \left(\frac{\lambda}{C}\right)^{\frac{2}{p+1}},
\end{eqnarray*}
for any $u\in H$. Therefore
\begin{eqnarray*}
M_{\lambda}\geq \left(\frac{\lambda}{C}\right)^{\frac{2}{p+1}}>0.
\end{eqnarray*}
$\Box$

Two additional properties will be used to prove the existence result:
\begin{itemize}
\item From (\ref{amns_e11}) (see also formula (17) in \cite{EsfahaniL2013}, along with (\ref{amns_e10}) of the present paper) we obtain the coercivity of $I(u)$, which is then equivalent to $||u||_{H}^{2}$ in all the cases (i) to (iv).
\item $M_{\lambda}$ is strictly subadditive in the sense that
\begin{eqnarray*}
M_{\lambda}\leq M_{\lambda_{1}}+M_{\lambda-\lambda_{1}},
\end{eqnarray*}
for any $\lambda_{1}\in (0,\lambda)$. Actually (see Note \ref{remark32b}), as in \cite{EsfahaniL2013} (see also \cite{LiuV2004} for the case $f(u)=u^{2}/2$) we have 
\begin{eqnarray}
M_{\lambda}=\lambda^{\frac{2}{p+1}}M_{1},\label{amns_e18}
\end{eqnarray}
for all $\lambda>0$.
\end{itemize}
The main result of existence is the following.
\begin{theorem}
\label{teorem1}
Under any of the conditions (i) to (iv) of Lemma \ref{lemma2}, let $\lambda>0$ and $\{u_{n}\}_{n}$ be a minimizing sequence in $H$ for $\lambda$. Then there exist subsequences $\{u_{n}\}_{n}$ in $H$, $\{y_{n}\}_{n}$ in $\mathbb{R}$ and $u\in H$ such that $u_{n}(\cdot+y_{n})\rightarrow u$ strongly in $H$. Furthermore, the function $u$ achieves the minimum $I(u)=M_{\lambda}$ subject to $K(u)=\lambda$.
\end{theorem}
\begin{note}
The proof is similar to that of other references, see in particular \cite{EsfahaniL2013,Levandosky1997,LiuV2004}. 
\end{note}

{\proof} From coercivity of $I$, the sequence $\{u_{n}\}$ is bounded in $H$ and therefore if 
%
we consider the $L^{1}$ sequence
\begin{eqnarray*}
\rho_{n}=|D_{x}^{1/2}u_{n}|^{2}+|\partial_{x}^{-1}u_{n}|^{2}+|\partial_{x}u_{n}|^{2}.
\end{eqnarray*}
then $\rho_{n}$ is bounded in $L^{1}$. Thus, there is a subsequence $\{\rho_{n}\}$ with 
\begin{eqnarray*}
L=\lim_{n\rightarrow\infty}||\rho_{n}||_{L^{1}},
\end{eqnarray*}
and normalizing (by taking $\tilde{\rho}_{n}(x)=L\rho_{n}\left(||\rho_{n}||_{L^{1}}x\right)$) we may assume $||\rho_{n}||_{L^{1}}=L$ for all $n$.

If we apply the Concentration-Compactness Lemma, \cite{Lions}, to $\rho_{n}$ we have three possibilities:
\begin{itemize}
\item[(a)] Compactness: there exist $y_{k}\in\mathbb{R}$ such that for any $\epsilon>0$ there is $R(\epsilon)>0$ such that for all $k$
\begin{eqnarray*}
\int_{|x-y_{k}|\leq R(\epsilon)}\rho_{k}dx\geq \int_{-\infty}^{\infty}\rho_{k}dx-\epsilon=L-\epsilon.
\end{eqnarray*}
\item[(b)] Vanishing: For every $R>0$
\begin{eqnarray*}
\lim_{k\rightarrow\infty}\sup_{y\in\mathbb{R}}\int_{|x-y|\leq R}\rho_{k}dx=0.
\end{eqnarray*}
\item[(c)] Dichotomy: there exists $l\in (0,L)$ such that for all $\epsilon>0$ there are $R,R_{k}\rightarrow\infty, y_{k}\in\mathbb{R}$ and $k_{0}$ satisfying
\begin{eqnarray*}
\left|\int_{|x-y_{k}|\leq R}\rho_{k}dx-l\right|<\epsilon,\quad 
\left|\int_{R<|x-y_{k}|\leq R_{k}}\rho_{k}dx\right|<\epsilon,
\end{eqnarray*}
for $k>k_{0}$.
\end{itemize}

The next step is ruling out possibilities (b) and (c). Here the arguments are similar to those of, for example, \cite{Levandosky1997}. Assume that (b) holds. Using (\ref{amns_e16}) and the homogeneity of $F$, we have
\begin{eqnarray*}
\int_{|x-y|\leq 1}F(u_{n})dx\leq C\int_{|x-y|\leq 1}|u_{n}|^{p+1}dx
\leq C\left(\int_{|x-y|\leq 1}\rho_{n}dx\right)^{\frac{p+1}{2}},
\end{eqnarray*}
for all $y\in \mathbb{R}$ and some constant $C$. By (b) one can choose $n(\epsilon)$ so large that
\begin{eqnarray*}
\int_{|x-y|\leq 1}F(u_{n})dx
\leq C\epsilon^{\frac{p-1}{2}}\int_{|x-y|\leq 1}\rho_{n}dx,
\end{eqnarray*}
for $n\geq n(\epsilon)$. Summing over intervals centered at even integers $y=2k$ we have
\begin{eqnarray*}
K(u_{n})\leq C\epsilon^{\frac{p-1}{2}},\quad n\geq n(\epsilon),
\end{eqnarray*}
that is $K(u_{n})\rightarrow 0$ as $n\rightarrow\infty$, which is in contradiction with the assumption that $u_{n}$ is a minimizing sequence.

Assume now that (c) holds. Define cutoff functions $\xi_{1},\xi_{2}$ with support on $|x|\leq 2$ and $|x|\geq 1/2$ respectively and with $\xi_{1}(x)=1, |x|\leq 1, \xi_{2}(x)=1, |x|\geq 1$. Let us consider
\begin{eqnarray*}
u_{k,1}(x)&=&\xi_{1}\left(|x-y_{k}|/R\right)u_{k}(x),\\
u_{k,2}(x)&=&\xi_{2}\left(|x-y_{k}|/R_{k}\right)u_{k}(x).
\end{eqnarray*}
Then $u_{k,j}, j=1,2$ satisfy, for $k\geq k_{0}$
\begin{eqnarray*}
I(u_{k})&=&I(u_{k,1})+I(u_{k,2})+O(\epsilon),\\
K(u_{k})&=&K(u_{k,1})+K(u_{k,2})+O(\epsilon).
\end{eqnarray*}
Since $u_{k}$ is bounded in $H$, then $u_{k,1}, u_{k,2}$ are bounded in $H$ independently of $\epsilon$. Therefore $K(u_{k,1}), K(u_{k,2})$ are bounded and then there are subsequences such that
\begin{eqnarray*}
\lambda_{i}(\epsilon)=\lim_{k\rightarrow\infty}K(u_{k,i}),\quad i=1,2,
\end{eqnarray*}
where $\lambda_{i}(\epsilon), i=1,2$ are bounded independently of $\epsilon$. So there is a sequence $\epsilon_{j}\rightarrow 0$ such that $\lambda_{i}(\epsilon_{j})\rightarrow \lambda_{i}$ for some $\lambda_{i}$ with $\lambda_{1}+\lambda_{2}= \lambda$. This leads to three possibilities:
\begin{enumerate}
\item $\lambda_{1}\in (0,\lambda)$. Then, by (\ref{amns_e18})
\begin{eqnarray*}
I(u_{k})&=&I(u_{k,1})+I(u_{k,2})+O(\epsilon_{j})
\geq M_{K(u_{k,1})}+M_{K(u_{k,2})}+O(\epsilon_{j})\\
&=&\left(K(u_{k,1})^{\frac{2}{p+1}}+K(u_{k,2})^{\frac{2}{p+1}}\right)M_{1}+O(\epsilon_{j})
\end{eqnarray*}
Taking $k\rightarrow\infty$, using that $u_{k}$ is a minimizing sequence and (\ref{amns_e18}), we have
\begin{eqnarray*}
\lambda^{\frac{2}{p+1}}M_{1}=M_{\lambda}\geq 
\left((\lambda_{1}(\epsilon_{j}))^{\frac{2}{p+1}}+(\lambda_{2}(\epsilon_{j}))^{\frac{2}{p+1}}\right)M_{1}+O(\epsilon_{j}).
\end{eqnarray*}
And, finally, if $j\rightarrow\infty$ then
\begin{eqnarray*}
M_{1}\geq 
\left(\left(\frac{\lambda_{1}}{\lambda}\right)^{\frac{2}{p+1}}+\left(\frac{\lambda_{2}}{\lambda}\right)^{\frac{2}{p+1}}\right)M_{1}>M_{1},
\end{eqnarray*}
which is a contradiction.
\item $\lambda_{1}=0$ (the same applies to $\lambda_{2}$ when $\lambda_{1}=\lambda$). From coercivity of $I$, the definition of $u_{k,1}$ and the hypothesis of dichotomy, we have, for some constant $C$
\begin{eqnarray*}
I(u_{k,1})&\geq &C\int_{-\infty}^{\infty}\left( |D_{x}^{1/2}u_{k,1}|^{2}+|\partial_{x}u_{k,1}|^{2}+|\partial_{x}^{-1}u_{k,1}|^{2}\right)dx\\
&=&C\int_{|x-y_{k}|\leq 2R}\left( |D_{x}^{1/2}u_{k}|^{2}+|\partial_{x}u_{k}|^{2}+|\partial_{x}^{-1}u_{k}|^{2}\right)dx\\
&=&C(l+O(\epsilon_{j})).
\end{eqnarray*}
Then 
\begin{eqnarray*}
I(u_{k})\geq C(l+O(\epsilon_{j}))+K(u_{k,2})^{\frac{2}{p+1}}M_{1}+O(\epsilon_{j}).
\end{eqnarray*}
As above, if $k\rightarrow\infty$ and then $j\rightarrow\infty$, we have
\begin{eqnarray*}
M_{1}\geq C
\left(\frac{l}{\lambda^{\frac{2}{p+1}}}\right)+\left(\frac{\lambda_{2}}{\lambda}\right)^{\frac{2}{p+1}}M_{1}>M_{1}.
\end{eqnarray*}
\item $\lambda_{1}>\lambda$ (the same applies to $\lambda_{2}$ if $\lambda_{1}<0$). Then, using the positivity of $I$ in all the cases (i) to (iv), we estimate
\begin{eqnarray*}
I(u_{k})\geq I(u_{k,1})+O(\epsilon_{j})\geq K(u_{k,1})^{\frac{2}{p+1}}M_{1}+O(\epsilon_{j}).
\end{eqnarray*}
And, again, if $k, j\rightarrow\infty$, then
\begin{eqnarray*}
M_{1}\geq\left(\frac{\lambda_{1}}{\lambda}\right)^{\frac{2}{p+1}}M_{1}>M_{1}.
\end{eqnarray*}
\end{enumerate}
So (c) is also ruled out and therefore compactness (a) holds. Now we prove that (a) implies the existence of a minimizer. (Again the arguments are similar to, e.~g. \cite{EsfahaniL2013,Levandosky1997}.) Since $u_{k}$ is bounded in $H$, there is a subsequence $u_{j}$ and $u\in H$ such that $\varphi_{j}=u_{j}(\cdot+y_{j})$ converges weakly to $u$ in $H$. Note also that, by Sobolev embedding, $u_{n}$ is bounded in $W^{1,q}(\mathbb{R})$ for all $q>2$, so the convergence is strong in $W^{1,p+1}_{loc}(\mathbb{R}), p>1$. Furthermore, by weak lower semicontinuity of $I$ in $H$, we have
\begin{eqnarray}
I(u)\leq \lim_{j\rightarrow\infty}I(\varphi_{j})=M_{\lambda}.\label{amns_e19}
\end{eqnarray}
Now we prove that $\varphi_{j}$ converges strongly to $u$ in $L^{p+1}$. We take $\sigma_{j}=|\varphi_{j}|^{p+1}$. From (\ref{amns_e18}) and compactness (a) of the $\rho_{k}$, $\sigma_{j}$ also satisfies (a). We take $\epsilon>0$ and $R_{0}>0$ so large that
\begin{eqnarray*}
\int_{|x|\geq R_{0}}|u|^{p+1}dx<\epsilon.\label{rmben20}
\end{eqnarray*}
By compactness of $\sigma_{j}$, there is $j_{1}(\epsilon)$ and $R(\epsilon)>R_{0}$ such that for $j>j_{1}(\epsilon)$
\begin{eqnarray}
&&\int_{|x|\geq R(\epsilon)}\sigma_{j}dx<\epsilon,\label{rmben21a}\\
&&\int_{|x|\geq R(\epsilon)}|u|^{p+1}dx\leq \int_{|x|\geq R_{0}}|u|^{p+1}dx<\epsilon.\label{rmben21b}
\end{eqnarray}
Therefore, if $B_{\epsilon}=B(0,R(\epsilon))$, by (\ref{rmben21a}), (\ref{rmben21b})
\begin{eqnarray*}
\int_{\mathbb{R}\backslash B_{\epsilon}}|\varphi_{j}-u|^{p+1}dx<2^{p+1}\epsilon,
\end{eqnarray*}
for $j>j_{1}(\epsilon)$. On the other hand, the strong convergence in $W^{1,p+1}_{loc}(\mathbb{R})$ implies the existence of $j_{2}(\epsilon)$ such that
\begin{eqnarray*}
\int_{ B_{\epsilon}}|\varphi_{j}-u|^{p+1}dx<\epsilon,
\end{eqnarray*}
for $j>j_{2}(\epsilon)$. Finally, if $j>\max\{j_{1},j_{2}\}$ then
\begin{eqnarray*}
\int_{\mathbb{R}}|\varphi_{j}-u|^{p+1}dx\leq
\int_{ B_{\epsilon}}|\varphi_{j}-u|^{p+1}dx+
\int_{\mathbb{R}\backslash B_{\epsilon}}|\varphi_{j}-u|^{p+1}dx
\leq (1+2^{p+1})\epsilon,
\end{eqnarray*}
which implies that $\varphi_{j}$ converges strongly to $u$ in $L^{p+1}$. Note that since $K$ is locally Lipschitz on $L^{p+1}$, \cite{Levandosky1997}, the strong convergence implies
\begin{eqnarray*}
K(u)=\lim_{j\rightarrow\infty}K(\varphi_{j})=\lambda.
\end{eqnarray*}
Therefore $I(u)\geq M_{\lambda}$ which, along with (\ref{amns_e19}), implies $I(u)=M_{\lambda}$, and $u$ is a minimizer of $I$ subject to $K(\cdot)=\lambda$. Finally, since $I$ is equivalent to $||\cdot ||_{H}^{2}$, $\varphi_{j}$ converges weakly to $u$ in $H$ and $I(\varphi_{j})\rightarrow I(u)=M_{\lambda}$, then $\varphi_{j}$ converges strongly to $u$ in $H$.$\Box$
\begin{note}
We observe that the conditions in Lemma \ref{lemma2} (or in Theorem \ref{teorem1}) seem to be in agreement with those of the limiting case $\beta=0$ (Ostrovsky equation) for the existence of solitary waves. More specifically, let us assume $\alpha=0$. If $\beta\rightarrow 0-$, then condition (i) in Lemma \ref{lemma2} implies $c_{s}<c^{*}=2\sqrt{\gamma\delta}$, while if $\beta\rightarrow 0+$ then
$z_{+}\rightarrow c^{*}$ and condition (iv) also reads $c_{s}<c^{*}$. This coincides with Theorem 2.1 of \cite{LevandoskyL2006} for the generalized Ostrovsky equation.
\end{note}
\begin{note}
A second observation is that the conditions in Lemma~\ref{lemma2} seem to be also in agreement with the arguments exposed n \cite{GalkinS1991,ObregonS1998} to justify the possibility of soliton solutions in the Ostrovsky equation, with $\alpha,\gamma,\delta>0$. Linearizing (\ref{rmben1}) and seeking for plane wave solutions, the corresponding dispersion relation for waves of small amplitude is 
\begin{equation}\label{amns_dispr}
\omega/k=\alpha+\frac{\gamma}{k^{2}}+\delta k^{2}-\beta |k|.
\end{equation}

\begin{figure}[h!]
\centering
{\includegraphics[width=6cm,height=6cm]{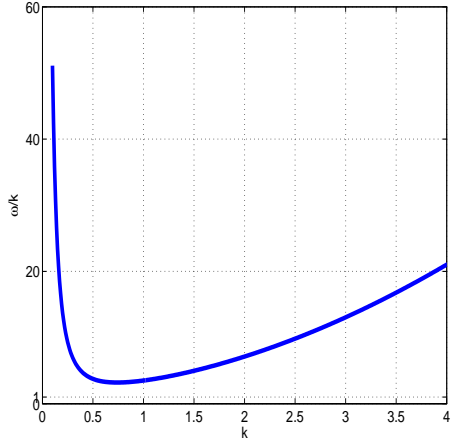}}
\caption{$\omega(k)/k$ vs $k$. Case $\beta<0$.}
\label{amns_fig7}
\end{figure}
The phase velocity (\ref{amns_dispr}) is displayed in Figure~\ref{amns_fig7} for $\beta<0$. in this case, linear perturbations can only exist within a semibounded range of phase velocities, giving the chance of having solitary waves which are not in resonance with linear perturbations and are not subject to radiative decay.

\begin{figure}[h!]
\centering
\subfigure[]
{\includegraphics[width=6cm,height=6cm]{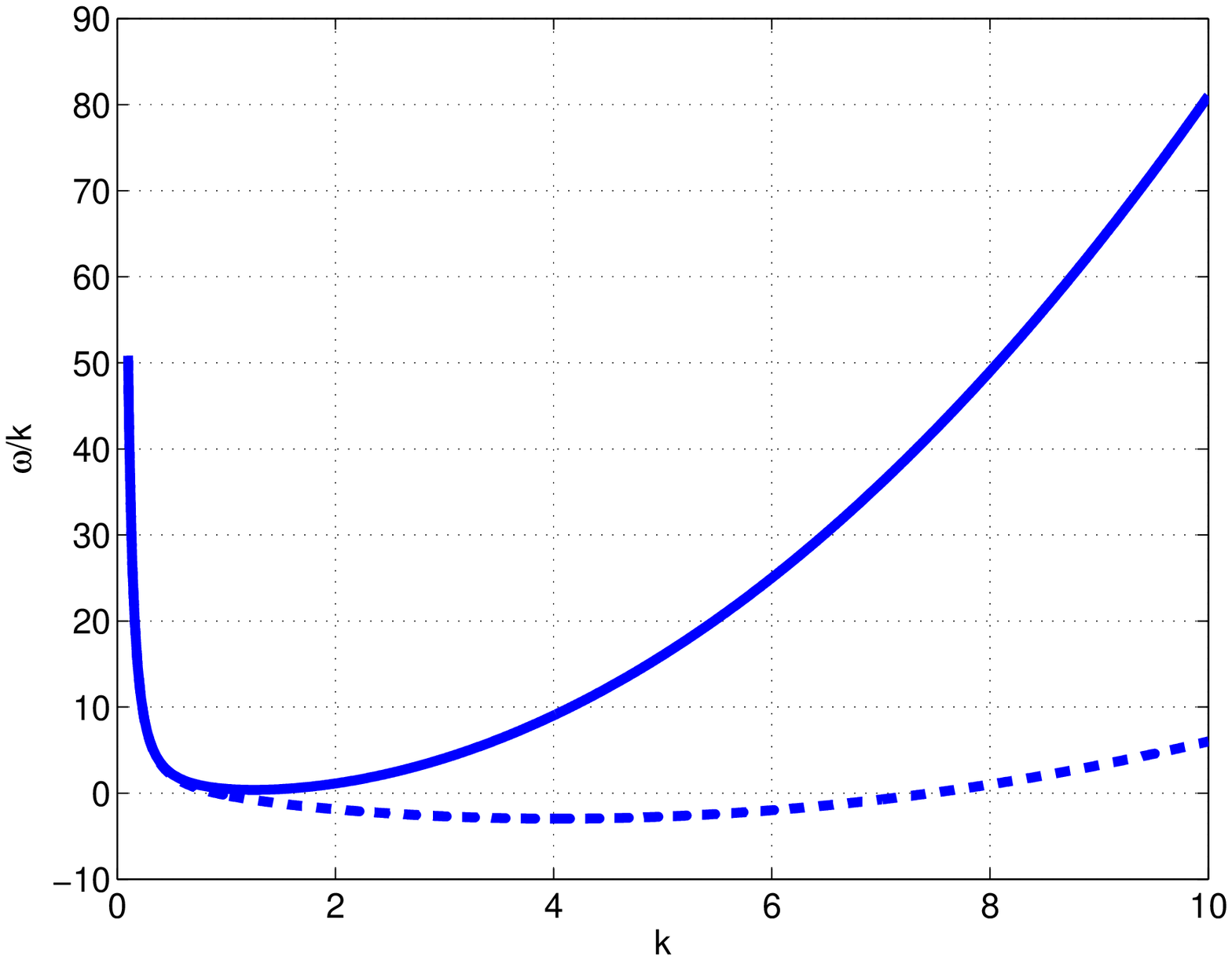}}
\subfigure[]
{\includegraphics[width=6cm,height=6cm]{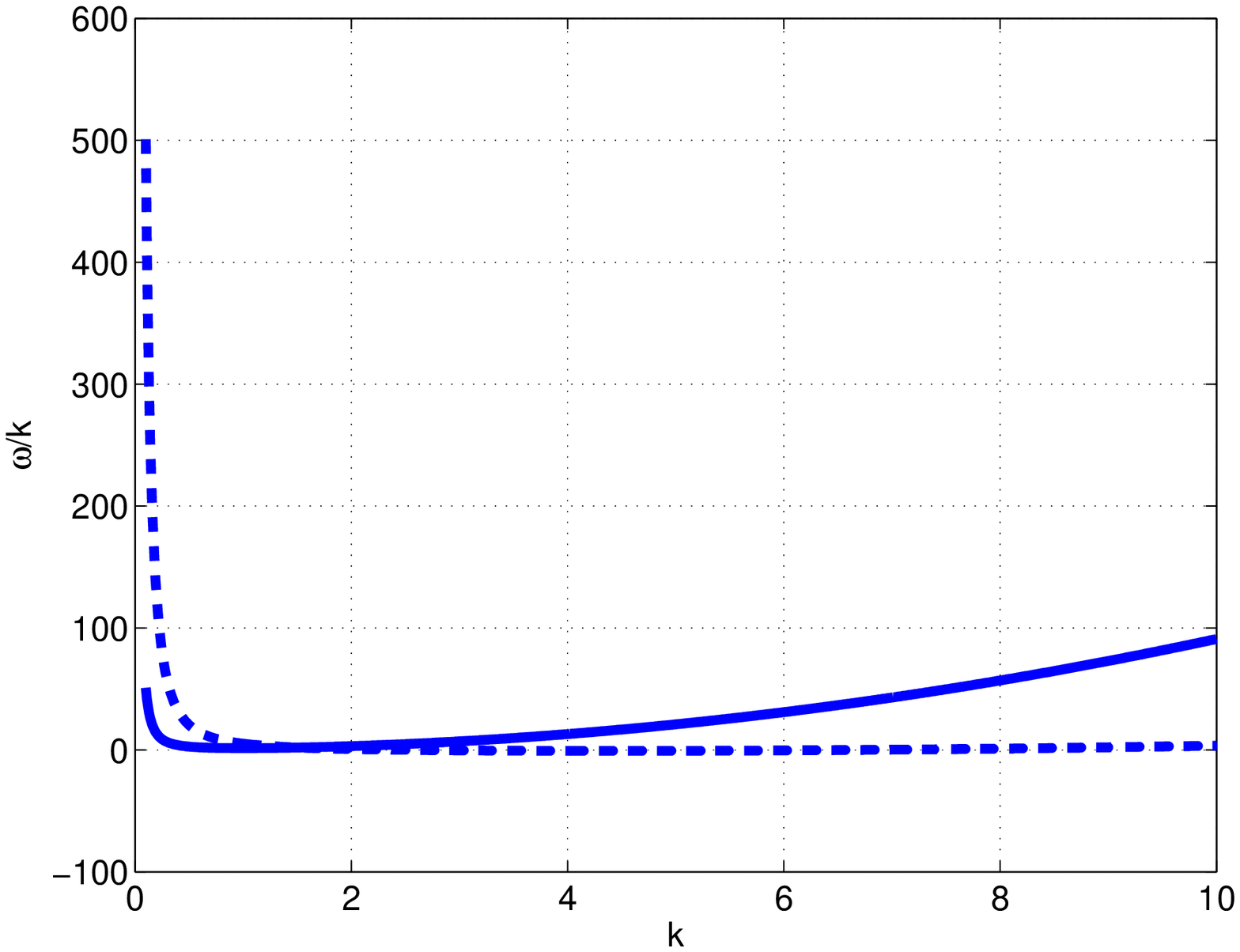}}
\subfigure[]
{\includegraphics[width=6cm,height=6cm]{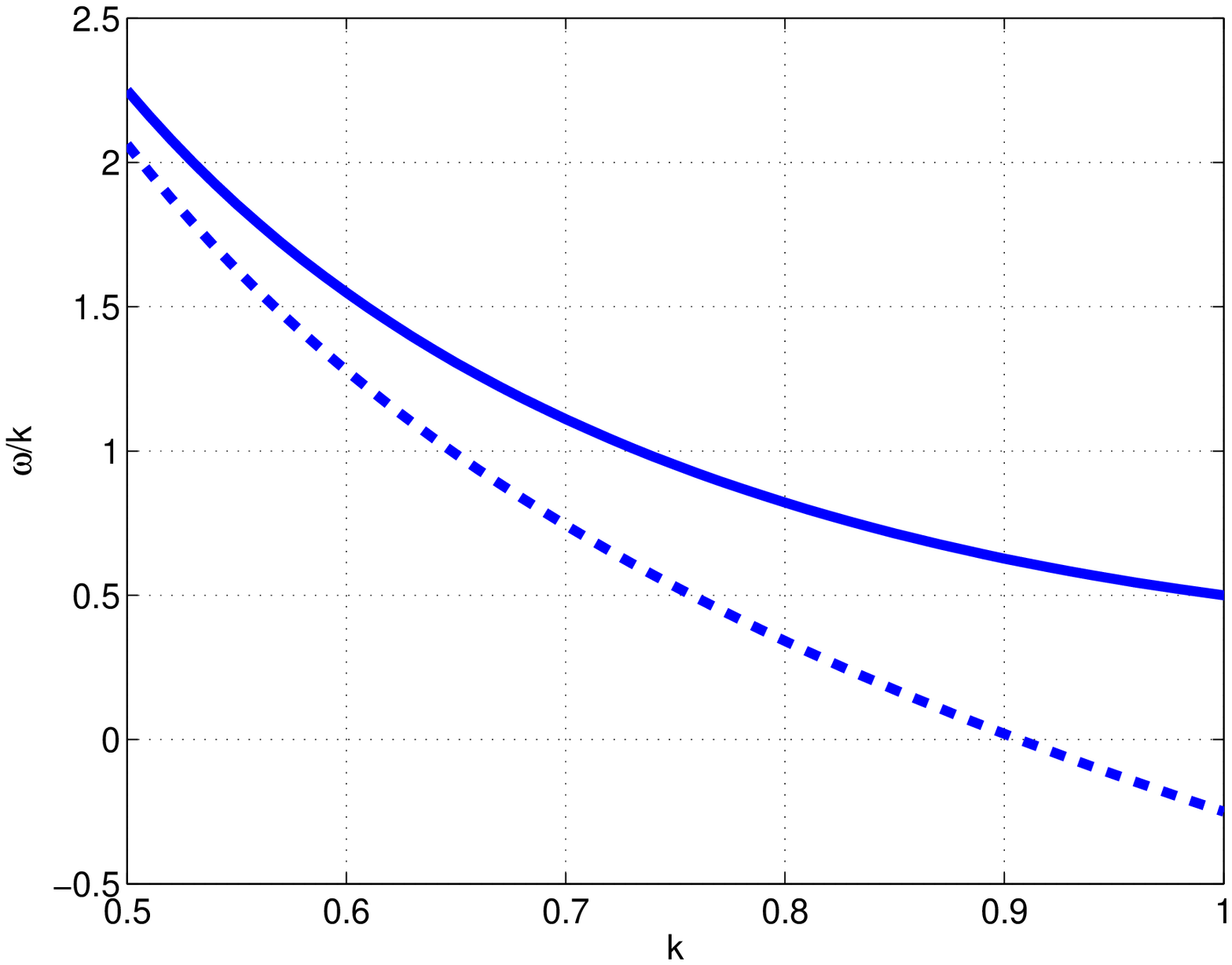}}
\subfigure[]
{\includegraphics[width=6cm,height=6cm]{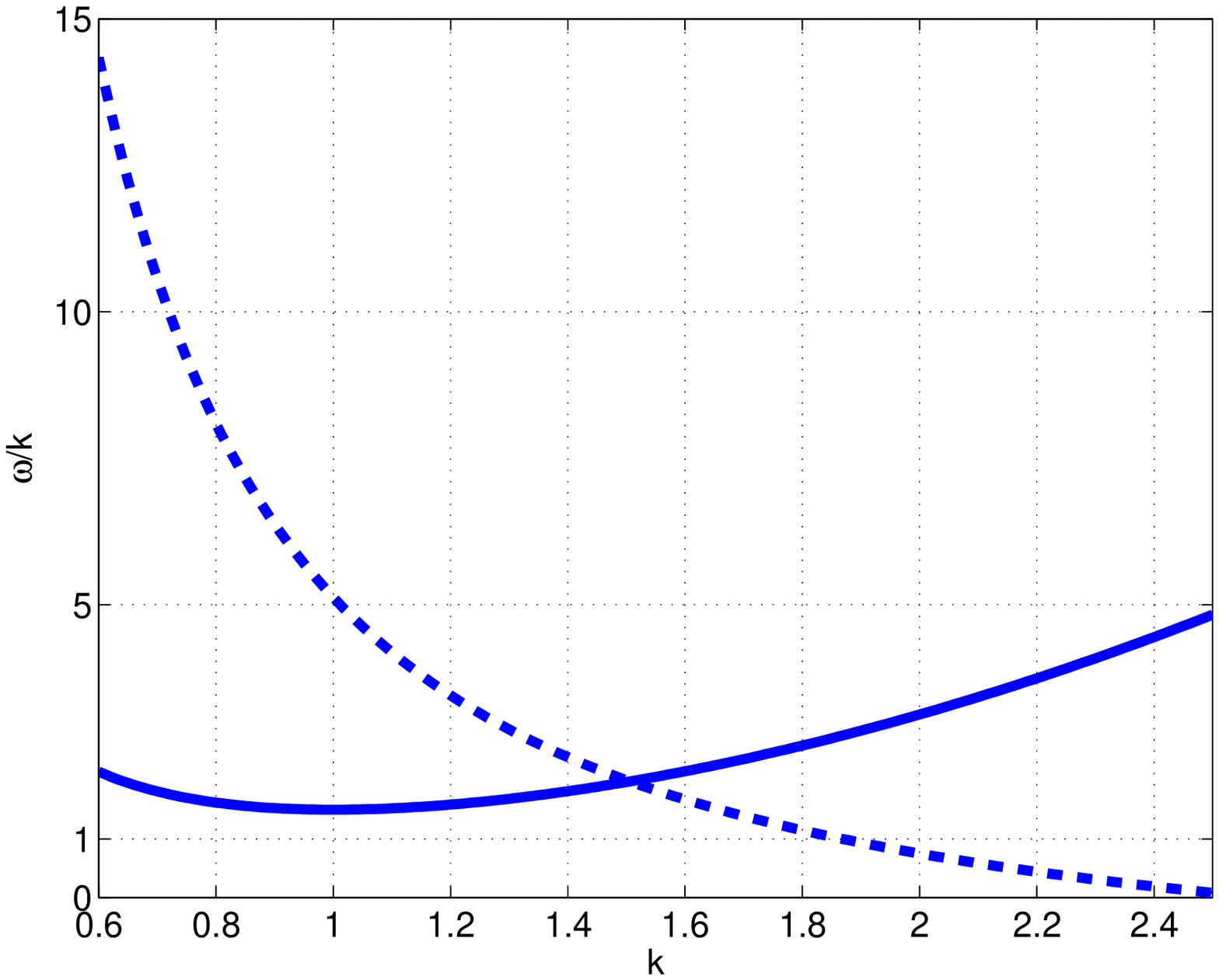}}
\caption{$\omega(k)/k$ vs $k$ for $\beta>0$. (a) Solid line: $A>0, B>0$ ($\alpha=\gamma=1/2, \beta=2,\delta=1$); dashed line: $A<0,B>0$ ($\alpha=\gamma=1/2, \beta=2,\delta=1/4$); (b) Solid line: $A>0, B<0$ ($\alpha=\gamma=1/2, \beta=\delta=1$); dashed line: $A<0,B<0$ ($\alpha=1/2,\gamma=5, \beta=1,\delta=1/8$); (c) Magnification of (a); (d) Magnification of (b)}
\label{amns_fig8}
\end{figure}
This relation, for $\beta>0$, is depicted in  Figure~\ref{amns_fig8},
depending on the sign of $A=4\delta-\beta>0$ and $B=\beta^{3}-\gamma(4\delta-\beta)^{2}$.
According to the arguments in \cite{GalkinS1991}, solitary waves are only possible when $A>0, B<0$ see solid lines in Figure~\ref{amns_fig8}(b) and (d).
\end{note}
\subsection{Comparisons with the Ostrovsky equation}
\begin{figure}[h!]
\centering
{\includegraphics[width=0.7\textwidth]{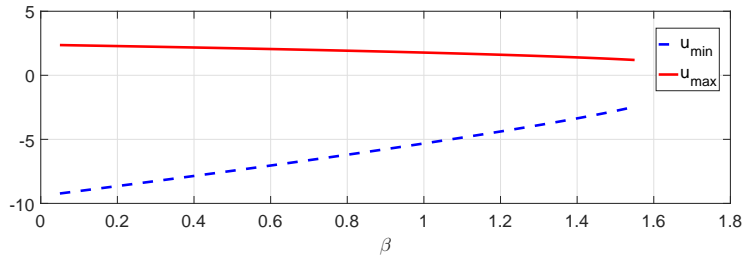}}
\caption{Amplitude vs $\beta$.}
\label{amns_fig9}
\end{figure}
The purpose of this section is to compare, by computational means and through the corresponding solitary waves, the equation (\ref{amns_e1}) with the Ostrovsky equation, the classical model for internal waves in rotating fluids, which is the limiting case of (\ref{amns_e1}) by taking $\beta=0$.

A first observation in this sense is concerned with the behaviour of the amplitude of the solitary waves of (\ref{amns_e1}) as function of $\beta$, illustrated in Figure~\ref{amns_fig9}. Note that the maximum positive excursion of the profiles $u_{max}$ is decreasing and the minimum negative excursion $u_{min}$ is increasing as $\beta$ grows. For fixed values of the rest of the parameters, the solitary wave solutions of the Ostrovsky equation ($\beta=0$) gives then the maximum amplitude.

\begin{figure}[h!]
\centering
{\includegraphics[width=0.7\textwidth]{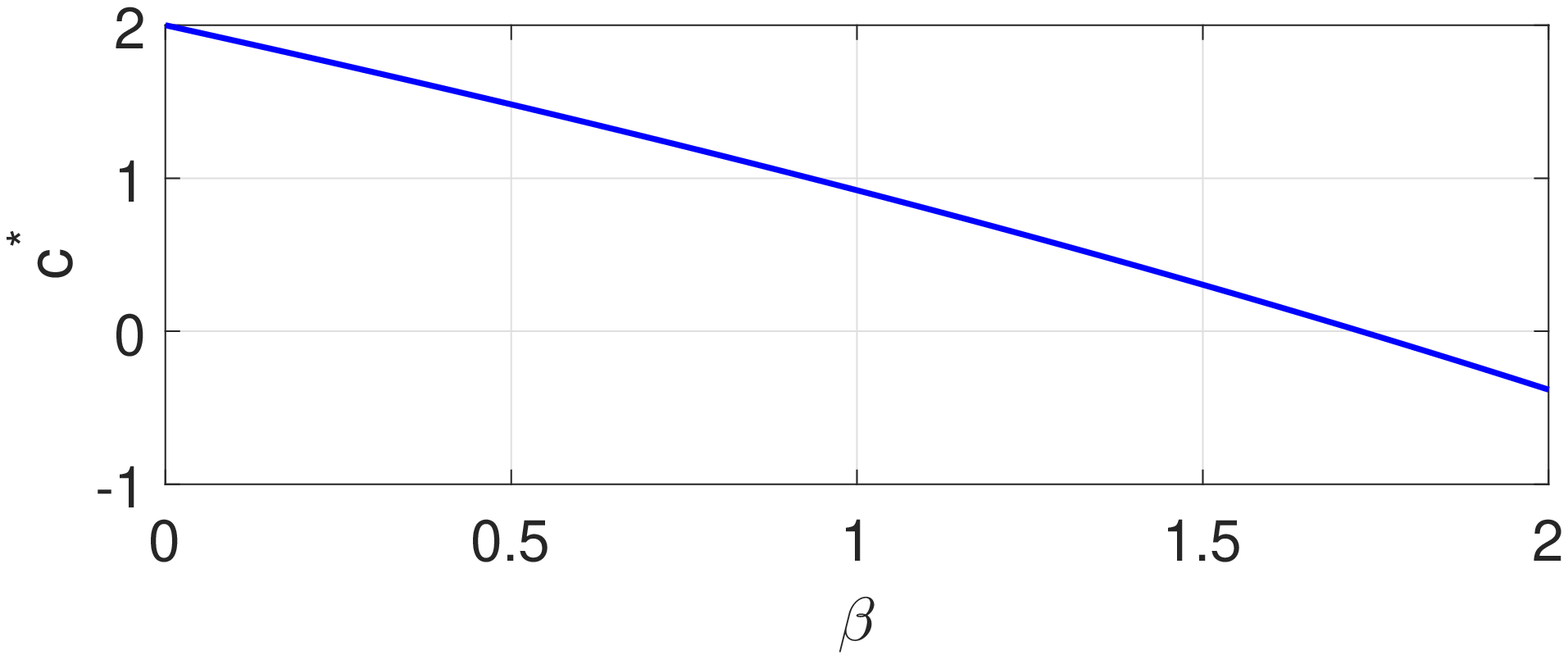}}
\caption{$c^{*}=\frac{1}{2}\left(-\beta\left(1+\frac{\beta}{4\delta}\right)+(4\delta-\beta)\sqrt{\frac{\gamma}{\delta}+\left(\frac{\beta}{4\delta}\right)^{2}}\right)$ vs $\beta$ with $\alpha=0, \gamma=\delta=1$.}
\label{amns_fig10}
\end{figure}

We also consder the behaviour of the limiting value of the speed to ensure the existence of solitary waves as function of $\beta>0$, predicted by Lemma~\ref{lemma2} and displayed in Figure~\ref{amns_fig10}. This shows that the maximum speed is a decreasing function of $\beta$ and then the maximum range of speeds to have solitary waves is given when $\beta=0$, that is, in the case of the Ostrovsky equation.

\begin{figure}[h!]
\centering
\subfigure[]
{\includegraphics[width=0.75\textwidth]{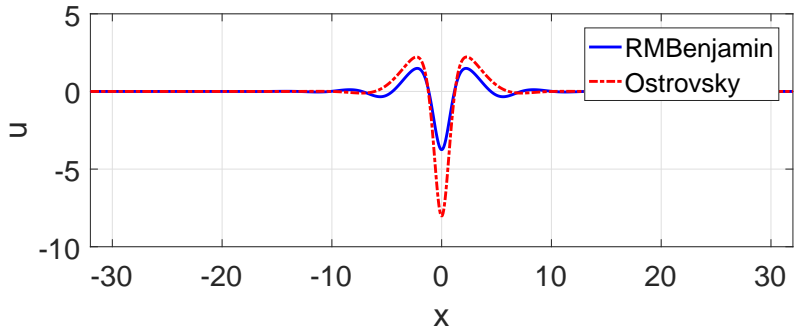}}
\subfigure[]
{\includegraphics[width=0.75\textwidth]{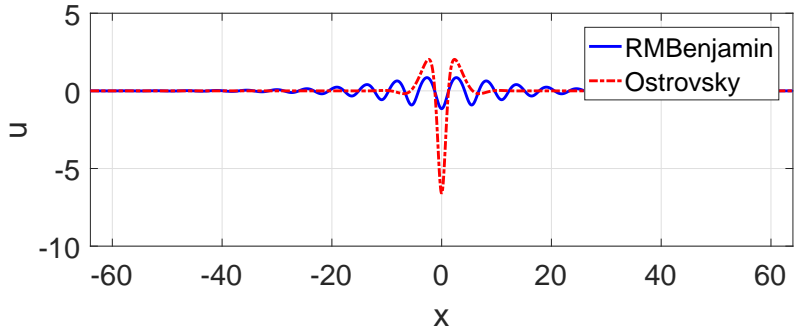}}
\caption{RMBenjamin vs Ostrovsky equations. Computed solitary wave profiles with (a) $c_{s}=0.1$, (b) $c_{s}=0.9$.}
\label{amns_fig11}
\end{figure}
\begin{figure}[h!]
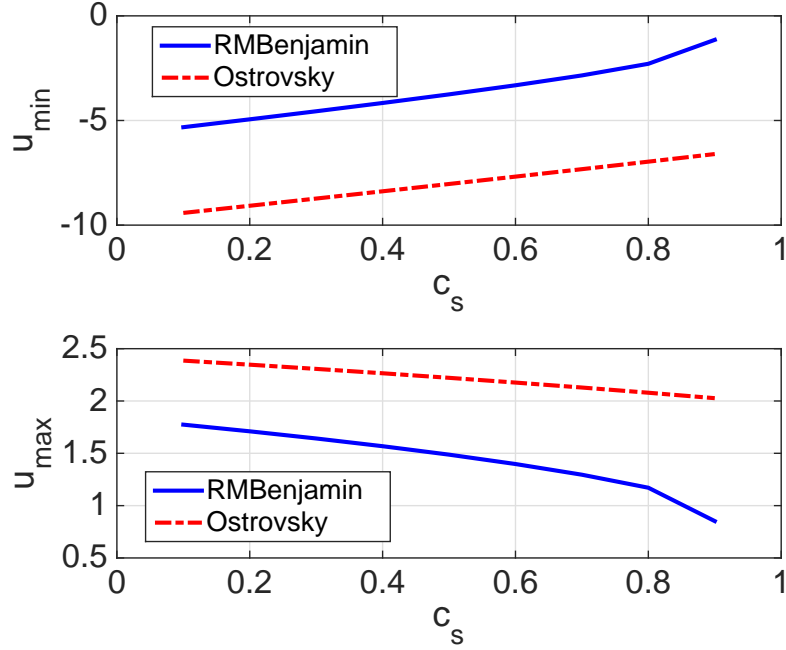

\centering
\subfigure
{\includegraphics[width=0.75\textwidth]{rmbenfig6b.eps}}
\subfigure
{\includegraphics[width=0.75\textwidth]{rmbenfig7b.eps}}
\caption{RMBenjamin vs Ostrovsky equations. Speed-amplitude relations.}
\label{amns_fig12}
\end{figure}

These two observations can explain the comparisons between the solitary waves of (\ref{amns_e1}) and of the Ostrovsky equation shown in Figures~\ref{amns_fig11} and ~\ref{amns_fig12}. Figure~\ref{amns_fig11} depicts the profiles corresponding to each equation for two values of the speed. According to this and the previous figures, note that the presence of the nonlocal term in (\ref{amns_e1}) with $\beta>0$ accelerates the formation of the oscillations in the profiles.

This is confirmed when we compare the behaviour of the maximum ($u_{max}$) and minimum ($u_{min}$) values of the waves as functions of the speed. Observe that in the case of the minimum, the value of the profile associated to the Ostrovsky equation is like a lower bound, while for the maximum, the corresponding value for the Ostrovsky equation is always above.
\section{Conclusions}
\label{sec4}
The present paper introduces a nonlinear dispersive nonlocal model for the propagation of internal waves in a two-layer system and under the presence of gravity, surface tension and rotational forces. The model can be derived from the inclusion of gravity effects in the rotating fluid model given by the Ostrovsky equation or by incorporating a dispersive, rotational component in the nonrotating model of the Benjamin equation. The proposed system can also be generalized  by including nonlinear terms from quadratic to any of homogeneous type with degree of homogeneity greater than two.

Three mathematical aspects of the model and its generalizations are analyzed. The first one is concerned with linear well-posedness and here sufficient conditions for existence and uniqueness of solution of the corresponsing IVP of the linear problem are established by using the theory developed in \cite{KenigPV1991b} and in terms of the parameters (with the corresponding physical meaning) of the equation. The second result is the derivation of three conservations laws and the Hamiltonian formulation, in accordance with its limiting cases of the Benjamin and Ostrovsky equations. Finally, the existence of solitary wave solutions is discussed, computationally and analytically. The generation of approximate solitary-wave profiles, described and developed in the present paper, gives a first indication of existence of solitary waves, suggests some of their properties (such as the amplitude-speed relation and the oscillatory decay) and allows to make comparisons between the proposed model and the classical rotating-fluid model given by the Ostrovsky equation. On the other hand, a theoretical result of existence, in terms of the parameters of the equation, is derived by using the Concentration-Compactness theory, \cite{Lions}, as in related rotating and nonrotating models.

Some open questions for a future research can be finally mentioned: 
\begin{itemize}
\item The first one is to make progress in the study of linear and nonlinear well-posedness, as well as in the proof of regularity and asymptotic decay of the solitary wave solutions, suggested by the numerical experiments.
\item Another research line is in the study of the stability of the solitary waves, both orbital and asymptotic, either theoretically or computationally. In this last point, the use of efficient numerical integrators, in order to have accurate long term simulations, is required.
\item A third open question is the analysis of the influence of the rotational effects from nonrotating-fluid models in more detail; in particular, it would be worth studying the weak rotation limit to the Benjamin equation, in a sort of comparison witrh the analogous property between the Ostrovsky equation and its weak rotation limit model, the KdV equation, \cite{LevandoskyL2007,LiuV2004,Tsuwaga2009}.
\end{itemize}

\section*{Acknowledgements}
This work was supported by FEDER and Junta de Castilla y Le\'on under the Grant VA041P17.


\begin{thebibliography}{999}
\bibitem{abr} J. P. Albert, J. L. Bona, J. M. Restrepo. (1999),
Solitary-wave solutions of the Benjamin equation, SIAM J. Appl.
Math., 59, pp. 2139-2161.
\bibitem{ApelOS2006} J. R. Apel, L. A. Ostrovsky, Y. A. Stepanyants. (2006), Internal solitons in the ocean, Technical Report.
\bibitem{ben0} T. B. Benjamin (1967), Internal waves of permanent form in
fluids of great depth, J. Fluid Mech. 29, pp. 559-592.
\bibitem{ben1} T. B. Benjamin. (1992), A new kind of solitary wave, J.
Fluid Mech., 245, pp. 401-411.
\bibitem{ben2} T. B. Benjamin. (1996), Solitary and periodic waves of a
new kind, Phil. Trans. R. Soc. Lond. A, 354, pp. 1775-1806.
\bibitem{Bona1980} J. L. Bona. (1980), Solitary waves and other phenomena associated with model equations for long waves. Fluid Dynamics Transactions 10 Panstowowe Wydawnictwo Naukowe: Warszawa, pp. 77-111.
\bibitem{Bona1981} J. L. Bona. (1981), Convergence of periodic wavetrains in the limit of large wavelength, Appl. Sci. Res. 37, pp. 21-30.
\bibitem{EsfahaniL2013} A. Esfahani, S. Levandosky. (2013), Solitary waves of the rotation-generalized Benjamin-Ono equation, Disc. Cont. Dyn. Sys., 33, pp. 663-700.
\bibitem{GalkinS1991} V. N. Galkin, Y. A. Stepanyants. (1991), On the existence of stationary solitary waves in a rotating fluid, J. Appl. Math. Mech., 55, pp. 939-943.
\bibitem{GilmanGS1995} O. A. Gilman, R. Grimshaw, Y. A. Stepanyants. (1995), Approximate analytical and numerical solutions of the stationary Ostrovsky equation, Stud. Appl. Math., 95, pp. 115-126.
\bibitem{GilmanGS1996} O. A. Gilman, R. Grimshaw, Y. A. Stepanyants. (1995), Dynamics of internal solitary waves in a rotating fluid, Dyn. Atm. Ocean, 23(1), pp. 403-411.
\bibitem{Grimshaw1985} R. H. Grimshaw. (1985), Evolution equations for weakly nonlinear, long internal waves in a rotating fluid, Stud. Appl. Math., 73, pp. 1-33.
\bibitem{Grimshaw1997} R. H. Grimshaw. (1997), Internal solitary waves, In P. L. F. Liu (ed.) Advances in Coastal and Ocean Engineering, World Scientific, Singapore, 3, pp. 1-30.
\bibitem{KenigPV1991b} C. E. Kenig, G. Ponce, and L. Vega. (1991), Oscillatory integrals and regularity of dispersive equations, Indiana Univ. Math. J. 40, pp. 33-69.
\bibitem{KenigPV1991} C.E. Kenig, G. Ponce, L. Vega. (1991), Well-posedness of the initial value problem for the Korteweg-de Vries equation, J. Amer. Math. Soc. 4, pp. 323-347.
\bibitem{KenigPV1993} C.E. Kenig, G. Ponce, L. Vega. (1993), Well-posedness and scattering results for the generalized Korteweg-de Vries equation via contraction principle, Comm. Pure Appl. Math. 46, pp. 527-620.
\bibitem{Levandosky1997} S.P. Levandosky. (1998), Stability and instability of fourth-order solitary waves, J. Dyn. Differential Equations, 10(1), pp. 151-188.
\bibitem{LevandoskyL2006} S. Levandosky, Y. Liu. (2006), Stability of solitary waves of a generalized Ostrovsky equation, SIAM J. Math. Anal., 38(3), pp.985-1011.
\bibitem{LevandoskyL2007} S. Levandosky, Y. Liu. (2007), Stability and weak Rotation limit of solitary
waves of the Ostrovsky equation, Disc. Cont. Dyn. Syst. Ser. B, 7(4), pp. 793-806.
\bibitem{LinaresM2004} F. Linares, A. Milan\'es. (2004), A note on solutions to a model for long internal waves in a rotating fluid, Mat. Contemp., 27, pp. 101-115.
\bibitem{LinaresM2006} F. Linares, A. Milan\'es. (2006), Local and global well-posedness for the Ostrovsky equation, J. Diff. Eq., 222, pp. 325-340.
\bibitem{Lions} P. L. Lions. (1984), The concentration-compactness principle in the calculus of variations. The locally compact case. Part I and Part II. Ann. Inst. Henri Poincar\'e Sect A (N.S.) 1, pp. 109-145 and pp. 223-283.
\bibitem{LiuV2004} Y. Liu, V. Varlamov. (2004), Stability of solitary waves and weak rotation limit for the Ostrovsky equation, J. Diff. Eq., 203, pp. 159-183.
\bibitem{ObregonS1998} M. A. Obregon, Y. A. Stepanyants. (1998), Oblique magneto-acoustic solitons in a rotating plasma, Phys. Lett. A, 249, pp. 315-323.
\bibitem{Ostrovsky1978} L. A. Ostrovsky. (1978), Nonlinear internal waves in a rotating ocean, Okeanologia, 18, pp. 181-191.
\bibitem{OstrovskyS1990} L. A. Ostrovsky, Y. A. Stepanyants. (1990), Nonlinear surface and internal waves in rotating fluids, in Nonlinear Waves 3, A. V. Gaponov-Grekhov, M. I. Rabinovich and J. Engelbrecht eds. Springer, New York, pp. 106-128.
\bibitem{PelinovskyS2004} {D. E. Pelinovsky and Y. A.
Stepanyants. (2004)}, { Convergence of Petviashvili's iteration method for
numerical approximation of stationary solutions of nonlinear wave
equations}, { SIAM J.~Numer.\ Anal.} { 42}, pp. 1110-1127.
\bibitem{Pet1976} { V. I. Petviashvili. (1976),} { Equation of
an extraordinary soliton}, { Soviet J.~Plasma Phys.} { 2}, pp.
257-258.
\bibitem{Shira1981} V. Shira. (1981), propagation of long nonlinear waves in a layer of a rotating fluid, Iza. Akad. Nauk SSSR, Fiz Atmosfery i Okeana, 17, pp. 76-81.
\bibitem{Shira1986} V. Shira. (1986), On long essentially nonlinear waves in a rotating ocean,  Iza. Akad. Nauk SSSR, Fiz Atmosfery i Okeana, 22, pp. 395-405.
\bibitem{Tsuwaga2009} K. Tsuwaga. (2009), Well-posedness and weak rotation limit for the Ostrovsky equation, J. Diff. Eq., 247, pp. 3163-3180.
\bibitem{VarlamovL2004} V. Varlamov, Y. Liu. (2004), Cauchy problem for the Ostrovsky equation, Disc. Dyn. Syst., 10, pp. 731-751.
\bibitem{yang2}  J. Yang. (2010), {\em Nonlinear Waves in Integrable and Nonintegrable Systems},
SIAM, Philadelphia.
 \end{thebibliography}
\end{document}